\documentclass[12pt]{amsart}
\usepackage{amsfonts,amssymb, amscd, latexsym, graphicx, psfrag, color, float}
\usepackage[all]{xy}

\newtheoremstyle{Herman}{}{}{}{}{}{.}{ }{}
\theoremstyle{Herman}
\newtheorem{po}{}
\newcommand \bp[1]{\begin{po}
#1\end{po}}

\begin{document}

\title[Cosimplicial $C^\infty$-rings and the de Rham complex of Euclidean space]{Cosimplicial $C^\infty$-rings and the de Rham complex of Euclidean space}

\begin{abstract}
A $C^\infty$-ring is a set equipped with $n$-ary operations corresponding to $n$-ary smooth functions $\mathbb{R}^n\rightarrow \mathbb{R}$. We prove that the cosimplicial abelian group associated to the de Rham complex of Euclidean space has the structure of a cosimplicial $C^\infty$-ring. We also analyse the notion of $R$-module (following Quillen) for a (co-)simplicial $C^\infty$-ring $R$.
\end{abstract}

\author{Herman Stel}
\address{Herman Stel, Dipartimento di Matematica e Informatica ``Ulisse Dini'',
 Viale Morgagni 67a, 50134 Firenze, Italy}
\email{herman\textunderscore stel@hotmail.com}

\maketitle

{\small \tableofcontents}

\section*{Introduction}

This note is based on work I did on $C^\infty$-rings as part of my Master's thesis \cite{St} at the University of Utrecht. 
It contains two theorems and one corollary. 
First, theorem \ref{DeRhamthm}, states that we have an isomorphism $N(\text{Inf}_n)\cong \text{dR}\left(\mathbb{R}^n\right)$ where $N$ is the normalization functor to cochain complexes, $\text{dR}\left(\mathbb{R}^n\right)$ is the de Rham complex of $\mathbb{R}^n$ and $\text{Inf}_n$ is a cosimplicial $C^\infty$-ring dual to the ``simplicial locus of infinitesimally small linear simplices in $R^n$''. I give explanation of this terminology in the corresponding section.
After this example we provide a second result, theorem \ref{mainsex}, stating that $T_{C^\infty\text{-ring}}\simeq C^\infty\text{-ring}\times_{\text{CRing}} T_{\text{CRing}}$, relating the tangent categories of the category of $C^\infty$-rings and of commutative rings. This is a recasting of proposition 2.2 in \cite{DK} in more modern language. It enables us to infer corollary \ref{impcor}, stating that modules in the sense of Quillen of a (co-)simplicial $C^\infty$-ring are cosimplicial modules of the underlying commutative ring. 

\bp{\textbf{Definition.} \textit{(Lawvere)} Write $\mathbb{E}$ for the full subcategory of the category of smooth manifolds spanned by the finite cartesian powers of $\mathbb{R}$, that is, spanned by $\mathbb{R}^n$ for $n\in \mathbb{Z}_{\geq 0}$. Then the category of $C^\infty$-rings is the full subcategory of the category of functors $\mathbb{E}\rightarrow \text{Set}$ spanned by the product preserving functors (called $C^\infty$-rings). $\Box$}

A $C^\infty$-ring $A$ is a set together with $n$-ary operations on that set corresponding to smooth functions on $\mathbb{R}^n$; the underlying set is $A(\mathbb{R})$. To get a feeling for this, note that every $C^\infty$-ring $A$ has an underlying commutative ring \[\left(A(\mathbb{R}),  A\left(+:\mathbb{R}^2\rightarrow \mathbb{R}\right), A\left(0:\{\ast\}\cong \mathbb{R}^0\rightarrow \mathbb{R}\right), A(-), A(\cdot),A(1)\right).\] The motivating example of a $C^\infty$-ring is the ring of smooth functions on a manifold $C^\infty\left(M,\mathbb{R}\right)$ and this fact gives a product preserving embedding $C^\infty\text{-Man}\rightarrow C^\infty\text{-ring}^{\text{op}}$. Recently $C^\infty$-rings have found a purpose in derived differential geometry. A good introduction to the theory of $C^\infty$-rings is the classical \cite{MR}. Dominic Joyce has incorporated an introduction in one of his books \cite{Jd-M} of which a preliminary copy is available from his website. The following theorem is pivotal to much of the theory of $C^\infty$-rings and some authors refer to it as Hadamard's Lemma.

\bp{\label{Had}\textbf{Taylor's Theorem.} Any infinitely differentiable function $f:U\rightarrow \mathbb{R}$ for $U$ a starlike neighborhood of a point $p\in \mathbb{R}^k$ satisfies, for any natural number $n$, the formula 
\[f(x)=\sum_{\tau\in T} \left(\prod_{i=1}^k\left( \frac{1}{\tau_i !}(x_i-p_i)^{\tau_i}\right) \frac{\partial^{\tau_1}}{\partial x_1^{\tau_1}}\cdots \frac{\partial^{\tau_k}}{\partial x_k^{\tau_k}}f\Big|_p\right) + \sum_{\sigma\in S}\left(\prod_{i=1}^{k}(x_i-p_i)^{\sigma_i}\right) g_\sigma(x),\]
for some $g_\sigma \in C^{\infty}(U,\mathbb{R})$; I wrote $T$ for $\{(\tau_1,...,\tau_k)\in \mathbb{Z}_{\geq 0}^k; \sum_{i=1}^k \tau_i \leq n\}$ and $S:=\{\sigma\in\mathbb{Z}_{\geq 0}^k;\sum_{i=1}^k \sigma_i=n+1\}$.\\
\underline{Proof:} Corollary 2.4 of \cite{N}. $\Box$}

It is important in our context that the $g_\sigma$ appearing in the formulation are smooth, and hence part of our algebraic theory of $C^\infty$-rings. One consequence is that for a $C^\infty$-ring $R$ the congruences on $R$ are in bijection with the ideals of the underlying commutative ring of $R$, which makes the theory that much more manageable. Another place where it is used is in lemma \ref{Coomods}. Exactly how determining this lemma is may be seen in the old \cite{DK} and the recent preprints \cite{CR1, CR2}.\\

Some words about the \textbf{relation to other work}. The first section finds its main inspiration in the work of Anders Kock \cite{K}. In that book Kock argues entirely synthetically using a theory that finds a natural interpretation inside a smooth topos. Such topoi are usually built from $C^\infty$-rings (see \cite{MR}). In Kock's book one finds the construction of the section on the De Rham complex synthetically. Here we give the same construction for Euclidean space, but we do it wholly inside the category of $C^\infty$-rings, since the logical tools used are only those whose semantics rely on the existence of limits. Technically, our result proves that when we start with Euclidean space, Kock's construction yields a simplicial representable object. Our objects, $C^\infty$-rings, have underlying abelian groups and so our context allows for the application of the Normalization functor: $\text{Ab}^{\Delta} \rightarrow \text{coChain}(\mathbb{Z})^+$; the result is isomorphic to the de Rham complex.\\

The second section describes the theory of modules over $C^\infty$-rings through the tangent category. Some of the main ideas that go into lemma \ref{Coomods} can be found in Proposition 2.2 of \cite{DK}. The further step of putting this in the context of the tangent category allows us to say things about modules over simplicial and cosimplicial $C^\infty$-rings.\\

As written above, this was work done as part of my MSc thesis at the University of Utrecht. It is a spin-off of Urs Schreiber's project \cite{Sc} (see section 4.5.1 ``$\infty$-Lie algebroids'' there). He suggested these topics to me and I \textbf{thank} him for that. Thanks go to Nicol\`o Sibilla and Urs Schreiber for helping me polish the final draft. Finally, I am grateful to the University of Florence for their support and flexibility.

\section*{Infinitesimal simplices and the de Rham complex}

There is a rudimentary theory of synthetic differential geometry in which one has as primitive a commutative ring $R$ that plays the role of the real line and one has equality and conjunction as logical connectives. This is a small part of what is used in \cite{K}, but it suffices for the construction in this paragraph. We show how to construct ``the simplicial locus of infinitesimal linear simplices of Euclidean $n$-space'', which is written $\ell\text{Inf}_n$. The rudimentary theory can be interpreted, not only in the toposes of \cite{MR}, but also in the category opposite to the category of $C^\infty$-rings, since the logical connectives used can be interpreted categorially there using finite limits. One then obtains a cosimplicial $C^\infty$-ring, which has an underlying cosimplicial abelian group. Taking the normalization of this cosimplicial abelian group one obtains the de Rham complex, this is the content of theorem \ref{DeRhamthm}.\\

I feel that when one employs some relatively exotic language, as I take categorial logic to be, one should be careful to make it easy for the reader to check that one obeys the rules. But to be completely formal is not doable. Therefore, I have chosen the following compromise. We work inside $C^\infty\text{-ring}^{\text{op}}$ and define subobjects there using its internal language; this is done in the numbers labeled `Synthetic'. Doing so will define a subobject: an isomorphism-equivalence class of monomorphisms. Each synthetic definition is accompanied by an algebraic statement that characterizes the domain of a representative of such a subobject and by an informal explication. This way I hope the work will be readable by non-experts too. Indeed, if one is interested in theorem \ref{DeRhamthm} and wishes to skip the whole synthetic part, one can do so simply by not reading it, and taking all propositions labeled `Algebraic' as definitions. The reader wishing to learn categorial logic is referred to \cite{JEl}.\\

I should mention again the book \cite{K} by Anders Kock. The synthetic construction below is taken from there.\\

In the following we pass from an object in $C^\infty$-ring to an object in $C^\infty\text{-ring}^\text{op}$ by writing $\ell$ in front of it. We can go back by putting $C^\infty$ in front of the object in $C^\infty\text{-ring}^\text{op}$. Formally nothing happens, but that is a coincidence.

\bp{\textbf{Definition.} $R$ is the object of $C^\infty\text{-ring}^{\text{op}}$ corresponding to the $C^\infty$-ring $\mathbb{R}^n\mapsto C^\infty\left(\mathbb{R},\mathbb{R}^n\right)$.}

Since we are dealing with both wedge products and conjunctions formally, I will reserve $\wedge$ for the product and $\&$ for the conjunction. So $p\& q$ stands for $p \text{ and } q$ and $\&_{i} p_i$ stands for the conjunction of all $p_i$. 

\bp{\textbf{Definition.} \textit{(Synthetic)} Define $D(n)$ as the subobject of $R^n$ that corresponds to the formula $\&_{k=1}^n \&_{l=1}^n x_k\cdot x_l = 0.$ $\Box$}

In words, $D(n)$ is the subobject of $n$-tuples of points in $R$ whose pairwise product equals zero. There are several concepts that share the name infinitesimal but I here stipulate that $D(n)$ is the subobject of (first order) infinitesimally small vectors in $R^n$. If we take a representative of this subobject and consider its domain, we obtain a $C^\infty$-ring. As promised, I give an algebraic description as well.

\bp{\textbf{Proposition.}\label{algtild} \textit{(Algebraic description)} The domain of any monomorphism representing $D(n)$ is isomorphic, as a $C^\infty$-ring, to $C^\infty\left(\mathbb{R}^n, \mathbb{R}\right)/I_n$, where $I_n$ is the ideal generated by the functions $\pi_i\cdot \pi_j$ with $i$ and $j$ ranging over $\{1,...,n\}$, where $\pi_i: \mathbb{R}^n \rightarrow \mathbb{R}$ is the $i$-th projection. $\Box$}

In the following definition I allow myself some notational freedom: given a subobject $S$ of some object $X$ corresponding to a formula $\phi$ in a context containing the variable $x:X$, I write $x\in S$ for $\phi(x)$. Using that we get the following definition.

\bp{\textbf{Definition.} \textit{(Synthetic)} $\widetilde{D}(m,n)$ is the subobject of $\left(R^n\right)^m$ given by the formula
$\&_{j=1}^m x_j\in D(n) \; \& \;\&_{k=1}^m \&_{l=1}^m \left(x_k-x_l\in D(n)\right).$ $\Box$}

Informally, $\widetilde{D}(m,n)$ is the subobject of $m$-tuples of infinitesimally small vectors in $R^n$ that are pairwise infinitesimally close together.

\bp{\textbf{Proposition.} \textit{(Algebraic description)} The domain of any monomorphism representing $\widetilde{D}(m,n)$ is isomorphic to \[\frac{C^\infty\left(\left(\mathbb{R}^n\right)^m, \mathbb{R}\right)}{\left(\sum_{j=1}^m I_n \circ \pi_j \right) + \left(\sum_{k,l\in \{1,...,m\}} I_n\circ \left(\pi_k-\pi_l\right)\right)},\] where $I_n$ is the same as in proposition \ref{algtild}. $\Box$}

The following is called the subobject of $m+1$-tuples of pairwise infinitesimally near vectors in $R^n$.

\bp{\textbf{Definition.} \textit{(Synthetic)} Write $R^n_{<m>}$ for the subobject of $\left(R^n\right)^{m+1}$ given by the formula
$\&_{k=1}^{m+1} \&_{l=1}^{m+1} x_k - x_l \in D(n).$ $\Box$}

\bp{\textbf{Proposition.} \textit{(Algebraic description)} The domain of any monomorphism representing $R^n_{<m>}$ is isomorphic to \[C^\infty((\mathbb{R}^n)^{m+1},\mathbb{R})/\sum_{k=1}^{m+1}\sum_{l=1}^{m+1} I_n\circ (\pi_k-\pi_l).\quad \Box\]}

We can relate $\widetilde{D}(m,n)$ to $R^n_{<m>}$ as follows.

\bp{\textbf{Proposition.}\label{compare} For $n,m\in\mathbb{Z}_{\geq 0}$ we have $R^n\times \widetilde{D}(m,n)\cong R^n_{<m>}$.\\
\underline{Proof:} Working out the definitions and using the notation for the subobjects to denote some domain of a monomorphism representing them we have
\[R^n\times \widetilde{D}(m,n)\cong \ell\left(C^\infty((\mathbb{R}^n)^{m+1},\mathbb{R})/\sum_{k=2}^{m+1} \sum_{l=2}^{m+1} I\circ (\pi_k-\pi_l)+ \sum_{k=2}^{m+1} I\circ \pi_k\right)\]
and 
\[R^n_{<m>}\cong \ell\left(C^\infty((\mathbb{R}^n)^{m+1},\mathbb{R})/\sum_{k=1}^{m+1}\sum_{l=1}^{m+1} I\circ (\pi_k-\pi_l)\right),\] where $I:=\langle \pi_i\cdot \pi_j;\; i,j=1,...,n\rangle\subseteq C^\infty(\mathbb{R}^n,\mathbb{R})$. The image under $\ell C^\infty(-,\mathbb{R})$ of $\phi:(\mathbb{R}^n)^{m+1}\rightarrow (\mathbb{R}^n)^{m+1}: (x_1,...,x_{m+1})\mapsto (x_1, x_2-x_1, ..., x_{m+1}-x_1)$ then induces an isomorphism from $R_{<m>}^n$ to $R^n\times \widetilde{D}(m,n)$. One sees this by checking that $C^\infty\left(\phi^{-1},\mathbb{R}\right)$ is a two-sided inverse of the dual map. $\Box$}

In our rudimentary internal language we cannot talk about a simplicial object $S$ in $C^\infty$-ring as a whole. What we can do is talk about the degeneration and face maps and the individual objects $S_n$ for every $n$ separately. Thus, we can define each face and degeneracy map synthetically on $\left(R^n\right)^{m+1}$ and prove that if the input is in $R^n_{<m>}$, then the output is in $R^n_{<m-1>}$ (in the case of a face map). I chose the following, algebraic, form for the definition because defining it synthetically would be hard to read and tedious to write down. In the example below I show what this looks like synthetically.

\bp{\textbf{Definition.} \textit{(Algebraic)} Write $\overline{n}:=\{1,...,n\}$ and $F:\text{Set}\rightarrow C^\infty\text{-ring}$ for the free $C^\infty\text{-ring}$ functor. Then $F(\Delta([0],-)\times \overline{n})$ is a cosimplicial $C^\infty\text{-ring}$ and \[\mathcal{I}:[m]\mapsto \sum_{k=1}^{m+1}\sum_{l=1}^{m+1} I_n\circ (\pi_k-\pi_l)\] is a cosimplicial ideal of $F(\Delta([0],-)\times \overline{n})$. Write $\text{Inf}_n$ for $F(\Delta([0],-)\times \overline{n})/\mathcal{I}$. $\Box$}

\bp{\textbf{Definition.} $\ell\text{Inf}_n$ is the functor $\Delta^{\text{op}}\rightarrow \left(C^\infty\right)^{op}$ dual to $\text{Inf}_n$. $\Box$}

\bp{\textbf{Proposition.} $\left(\ell\text{Inf}_n\right)_m\rightarrow \left(R^n\right)^{m+1}$ is a representative for $R_{<m>}^n$.\\
\underline{Proof:} Note that $F\left(\Delta([0],[m]) \times \overline{n}\right)$ is isomorphic to $C^\infty\left(\left(\mathbb{R}^{n}\right)^{m+1}, \mathbb{R}\right)$, see \cite{MR}. Now one compares to the algebraic description of $R^n_{<m>}$. $\Box$}

To get the flavour of what this simplicial locus looks like, I give some examples of provable sequents.

\bp{\textbf{Example.} \textit{(Synthetic)} The following sequents are provable using the theory of finite limits.
\[\left(x_0,...,x_m\right)\in R^n_{<m>}\vdash_{x_0,...,x_m:R^n} \left(x_0,...,\widehat{x}_i,...,x_m\right)\in R^n_{<m-1>}\]
\[\left(x_0,...,x_m\right)\in R^n_{<m>}\vdash_{x_0,...,x_m:R^n} \left(x_0,...,x_i,x_i,...,x_m\right)\in R^n_{<m+1>}\quad \Box\]}

Classically we have a simplicial manifold of linear simplices $[m]\mapsto \left(\mathbb{R}^n\right)^{m+1}$ and this object is also present in our category. The subobjects $R^n_{<m>}$ we have defined  should therefore be viewed as an object of linear simplices, and we say that they are the infinitesimal linear $m$-simplices of $R^n$. The theorem mentions the normalization of cochain-complexes of abelian groups. The construction is dual to the one on page 218 of \cite{DP}.

\bp{\textbf{Theorem.} \label{DeRhamthm} Let $n\in \mathbb{Z}_{\geq 0}$ be a natural number. The normalized cochain complex of the underlying cosimplicial abelian group of $\text{Inf}_n$ is isomorphic to the De Rham complex of $\mathbb{R}^n$.\\
\underline{Proof:} I produce an isomorphism \[C^\infty\left(R^n\times \widetilde{D}(m,n)\right)\cong C^\infty\text{Man}\left(\mathbb{R}^n, \bigwedge_{i=1}^m \left(\mathbb{R}^n\right)^\ast\right)=:A_m,\] where $\bigwedge_{i=1}^m \left(\mathbb{R}^n\right)^\ast$ denotes the space of alternating tensors. The latter is isomorphic to the $\mathbb{R}$-vector space $\Omega^m\left(\mathbb{R}^n\right)$ of De Rham $m$-forms on $\mathbb{R}^n$, the former to $\text{Inf}_n$ by proposition \ref{compare}.\\
To compute the normalized cochain complex of $\text{Inf}_n$ we will first analyse the cofaces, some of whose images we will divide out by. The cosimplicial structure was obtained from the one on $C^\infty R^n_{<->}$ by transport of structure. One checks that for $i\geq 1$ the $(i\text{-th})$ coface map of the corresponding cosimplicial $C^\infty\text{-ring}$ sends $\pi_k$ to $\pi_k$ if $k<i$ and $\pi_k$ to $\pi_{k+1}$ if $k\geq i$. The $0\text{-th}$ coface map sends $\pi_1$ to $\pi_1+\pi_2$ and $\pi_k$ to $\pi_{k+1}$ for all $k>1$. We compute the normalized cochain complex of this cosimplicial $\mathbb{R}$-module by dividing out $C^\infty(R^n\times \widetilde{D}(m,n))$ by the sub-vectorspace generated by the images of all but the $0$-th coface map. To facilitate this procedure later on, we note that there is an equality of ideals
\[\begin{array}{llll}\sum_{l=2}^{m+1}\sum_{k=2}^{m+1} I\circ (\pi_k-\pi_l)+\sum_{k=2}^{m+1} I\circ \pi_k=\\\\ \langle \pi_{i,j}\pi_{i', j'} +\pi_{i,j'} \pi_{i', j};\; i,i'=2,...,m+1\; \& \; j,j'=1,...,n\rangle \end{array}\] as the reader can check for herself.\\
Now consider, for $[m]\in\text{Ob}\left(\Delta\right)$ the function
$\phi_{[m]}:C^\infty((\mathbb{R}^n)^{m+1},\mathbb{R})\rightarrow A_m$ that sends $f$ to $x\mapsto$
\[\sum_{\alpha:\{2,...,m+1\}\hookrightarrow \{1,...,n\}} \frac{\partial^\alpha f}{\partial \xi^\alpha}\Big|_{(x,0,...,0)} \cdot \bigwedge_{i=1}^{m} \pi_{\alpha(i+1)},\]   where the $\alpha:\{2,...,m+1\}\hookrightarrow \{1,...,n\}$ range over injections. Firstly, this is a well defined $\mathbb{R}$-linear function. Second, the ideal $\langle \pi_{i,j}\pi_{i', j'} +\pi_{i,j'} \pi_{i', j};\; i,i'=2,...,m+1\; \& \; j,j'=1,...,n\rangle$ is in the kernel of $\phi$ (we omit the subscript $[m]$). This one proves by using the definition of wedge product and noting that
\[x\mapsto \frac{\partial^\alpha (f\cdot ( \pi_{i,j}\pi_{i', j'} +\pi_{i,j'} \pi_{i', j}))}{\partial \xi^\alpha}\Big|_{(x,0,...,0)}\] is fixed under transposition of two of the non-$x$ variables. Also, the images of all but the $0$-th coface map are contained in ker$(\phi)$. Indeed, if some function $f$ does not depend on the ($1<$) $i$-th (vector valued) variable, then 
\[\frac{\partial^\alpha f}{\partial \xi^\alpha}\Big|_{(x,0,...,0)}\] will be zero. So this induces a linear map $\overline{\phi}$ from $N\left(C^\infty\left(R^n\times \widetilde{D}(-,n)\right)\right)_m$ to $A_m$. To go back, take a form $\omega$ and send it to the coset of $\psi(\omega):(x,y_1,...,y_m)\mapsto \omega(x)(y_1,...,y_m)$.\\
We must check that $\psi(\phi(f))$ is in the same coset as $f$. Taking the Taylor expansion of $(y_1,...,y_m)\mapsto f(x,y_1,...,y_m)$ up to order $m+1$ around $0$ and using Hadamard's lemma for the rest term we note that since $g\cdot \pi_{i,j}\pi_{i,j'}$ is in $\text{ker}(\phi)$ for smooth $g$ and since all terms of lower degree depend on less than $m$ variables of $y_1,...,y_m$, the function $f$ is equivalent modulo $\text{ker}\phi$, to
\[(x,y_1,...,y_m)\mapsto \sum_{\alpha:\{1,...,m\}\hookrightarrow \{1,...,n\}} \frac{\partial^\alpha f}{\partial \xi^\alpha}\Big|_{(x,0,...,0)}\cdot \prod_{i=1}^m y_{i,\alpha(i)}.\] The symmetrization of this map equals $\psi(\phi(f))$. Since for any smooth $g$ we have $g\cdot (\pi_{i,j}\pi_{i',j'}+\pi_{i,j'}\pi_{i',j})\in\text{ker}(\phi)$ we obtain $\psi(\phi(f))\sim f$.\\
Now we verify that $\phi(\psi(\omega))=\omega$ for any alternating $m$-form $\omega$. If $\omega$ is such a form then \[\frac{\partial^\alpha\psi(\omega)}{\partial\xi^\alpha}\Big|_{(x,0,...,0)}=\omega(x)(e_{1,\alpha(1)},...,e_{m,\alpha(m)}),\] where $e_{l,k}$ is a standard basisvector of $(\mathbb{R}^n)^m$; this can be checked expressing $(y_1,...,y_m)\mapsto\omega(x)(y_1,...,y_m)$ as a linear combination of maps of the form $\prod_{i=1}^m \pi_{i,p(i)}$ with $p:\{1,...,m\}\hookrightarrow \{1,...,n\}.$ Thus $\omega$ is equal to $x\mapsto$
\[(y_1,...,y_m)\mapsto \sum_{\alpha} \frac{\partial^\alpha\psi(\omega)}{\partial \xi^\alpha}\Big|_{(x,0,...,0)}\prod_{i=1}^m y_{i,\alpha(i)}.\] Its symmetrization is $\phi(\psi(\omega))$, but since each $\omega(x)$ is an alternating multilinear map $\phi(\psi(\omega))=\omega$.\\
I now show that $\phi\circ d^0\circ \psi=d$, where $d^0$ is the $0$-th face map and $d$ is the exterior derivative. Let $\omega\in A_m$. Then for some $\{ f_\alpha\in C^\infty(\mathbb{R}^n,\mathbb{R});\; \alpha:\{1,...,m\}\hookrightarrow \{1,...,n\}\}$ we have
\[\omega: x\mapsto \sum_{\alpha}\left(f_\alpha(x)\cdot \bigwedge_{i=1}^m \pi_{i,\alpha(i)}\right).\]
Thus, $d^0(\psi(\omega))=\sum_\alpha \left(f_\alpha \circ (\pi_1+\pi_2)) \cdot \bigwedge_{i=1}^m \pi_{i+2,\alpha(i)}\right)$
\[\sim_{\text{ker}(\phi_{m+1})} \sum_{\alpha}\left((f_\alpha\circ (\pi_1+\pi_2))\cdot \prod_{i=1}^m \pi_{i+2,\alpha(i)}\right)=:\chi.\]
We have $\phi_{m+1}(\chi)=\phi_{m+1}(d^0(\psi_m(\omega)))$. One computes, for $\beta:\{2,...,m+2\}\hookrightarrow \{1,...,n\}$, that 
\[\frac{\partial^\beta \chi}{\partial \xi^\beta}\Big|_{(x,0,...,0)}=\frac{\partial f_{\beta'}}{\partial \xi_{\beta(2)}}\Big|_{x},\]
where $\beta':\{1,...,m\}\hookrightarrow \{1,...,n\}:j\mapsto \beta(j+2)$. Consequently,
\[\phi_{m+1}(\chi):x\mapsto \sum_{\beta:\{2,...,m+2\}\hookrightarrow \{1,...,n\}} \frac{\partial f_{\beta'}}{\partial \xi_{\beta(2)}}\Big|_x \cdot \bigwedge_{i=1}^{m+1} \pi_{i,\beta(i+1)}=\]
\[\sum_{\alpha:\{1,...,m\}\hookrightarrow \{1,...,n\}} D_x f \wedge\bigwedge_{i=2}^{m+1} \pi_{i,\alpha(i-1)}=d(\omega)(x).\quad \Box\]}

\section*{Modules over $C^\infty$-rings}

This section starts from Quillen's definition of $R$-module for some object $R$ in a category with pullbacks. When applied to the category of commutative rings we obtain modules in the usual sense; we prove this well-known fact. We then apply said definition to the category of $C^\infty$-rings and show that what one obtains over a $C^\infty$-ring $R$ is the same as modules over the underlying commutative ring. In fact, we do everything somewhat more functorially: we consider the category of Quillen modules over all objects in the category, the tangent category. The main point of this section is theorem \ref{mainsex}; it has an important corollary \ref{impcor}.

\bp{\textbf{Definition.} Let $\mathcal{C}$ be a category with finite pullbacks. For an object $A\in \text{Ob}\left(\mathcal{C}\right)$ write $\text{Ab}(\mathcal{C}/A)$ for the category of abelian group objects in the category over $A$. We define the tangent category $T_\mathcal{C}$; an object is an object in $\text{Ab}(\mathcal{C}/A)$ for some $A$. If $X\in \text{Ab}(\mathcal{C}/A)$ and $Y\in \text{Ab}(\mathcal{C}/B)$ a morphism $X\rightarrow Y$ is a pair $(f,\phi)$ with $f:A\rightarrow B$ an arrow in $\mathcal{C}$ and $\phi:X\rightarrow f^\ast Y$ an arrow in $\text{Ab}(\mathcal{C}/A)$. $\Box$}

Urs Schreiber suggested the name `tangent category', which is inspired by Lurie's notion in \cite{L}
We show what one obtains in the case when $\mathcal{C}$ is the category $\text{CRing}$ of commutative rings.

\bp{\textbf{Definition.} We define a category $\text{Mod}$ of modules over commutative rings with multiplicative unit. As objects, take all modules over such rings. If $M$ is an $R$-module and $N$ is an $S$-module then a morphism $M\rightarrow N$ is a pair $(f,\phi)$ of a ring homomorphism $f:R\rightarrow S$ and an $R$-linear map $\phi:M\rightarrow N_f$, where $N_f$ is the $R$-module with structure map $R\rightarrow S\rightarrow \text{End}_{\text{Ab}}(N)$. $\Box$}

\bp{\textbf{Proposition.}\label{QuillenModuleRings} \textit{(Quillen)} The category of modules over commutative rings is equivalent to the tangent category $T_\text{CRing}$ of the category of commutative rings.\\
\underline{Proof:} We construct a functor $F:\text{Mod}\rightarrow T_\text{CRing}$ and show that it is full, faithful and surjective on objects.\\
Let us define the action of $F$ on objects. If $M$ is an $R$-module then $R\times M$ is an abelian group that can be given the structure of a ring by $(r_0,m_0)(r_1,m_1)=(r_0r_1, r_0m_1+r_1m_0)$. The projection $\pi_R:R\times M\rightarrow R$ can be equiped with an abelian group structure in $\text{Ring}/R$ by pulling back the group structure in Set of $M$ along $R\rightarrow 1$; one needs only check that this structure is in fact in $\text{Ab}(\text{Ring}/R)$, i.e. that the pullback of the arrows in Set become arrows in Ring. This defines $F$ on objects. If $M$ is an $R$-module and $N$ an $S$-module and $(f,\phi): M\rightarrow N$ is an arrow in $\text{Mod}$ then take $F(f,\phi)$ to be $(f,\phi')$, where $\phi': F(M)\rightarrow f^\ast(F(N)): (r,m)\mapsto (r,\phi(m))$.\\
It is clear that $F$ is full and faithful. Indeed, for $M$ an $R$-module and $N$ an $S$-module the inverse of $F_{M,N}$ is given by sending $(f,\phi')\mapsto (f,\pi\circ \phi')$, where $\pi: R\times_S (S\times N)\rightarrow S\times N\rightarrow N$ is the projection from the vertex of the pullback $f^\ast F(N)$ followed by the projection $S\times N\rightarrow N$.\\
To show that it is essentially surjective on objects, let $a:G\rightarrow R$ be an abelian group object in $\text{Ring}/R$ (by abuse of notation). Then the kernel of $a$ is an $R$-module and $\phi:G\rightarrow F(\text{ker}(a)):g\mapsto (a (g), g- \eta(a(g)))$ is an isomorphism of rings over $R$, where I wrote $\eta$ for the `unit element' of $a$. It is easily seen that it respects the group structure. $\Box$}

It turns out that when we apply Quillen's approach to modules to the category of $C^\infty$-rings we obtain the same thing, this is exactly expressed as follows. Some of the main ingredients of the following lemma can be found in \cite{DK}.

\bp{\textbf{Lemma.} \label{Coomods}Let $R$ be a $C^\infty$-ring and write $U(R)$ for its underlying ring. Then the forgetful functor $\mathcal{U}:\text{Ab}(C^\infty\text{-ring}/R)\rightarrow \text{Ab}(\text{Ring}/U(R))$ is an equivalence of categories.\\
\underline{Proof:} The hardest part of this proof is to show that the functor $\mathcal{U}$ is injective on objects. In the process of doing so we obtain a closed formula for the action of an abelian group in the overcategory of some $R$ on some smooth function $f:\mathbb{R}^k\rightarrow \mathbb{R}$ in terms of $R$ and the ring structure on the abelian group.\\
Let $R$ be a $C^\infty$-ring and let $a:M\rightarrow R$ be an abelian group in $C^\infty\text{-ring}/R$. Up to isomorphism $M$ sends $\mathbb{R}^n$ to $M(\mathbb{R})^n$. To see what it does on arrows, note that the underlying abelian group over the ring $U(R)$ is in the essential image of the functor $\text{Mod}\rightarrow T_{\text{Ring}}$ and the ring that is the domain of $\mathcal{U}(a)$ is therefore isomorphic to $R\times M_0$ (for some module $M_0$) with the ring structure as defined in the proof of lemma \ref{QuillenModuleRings}. Along this isomorphism we could write $(r,m)$ for a general element of $M(\mathbb{R})$, but I think it is better to save our formulas from too many commas and parentheses and write $r\oplus m$ for the same element.\\
Now let $f:\mathbb{R}^k\rightarrow \mathbb{R}$ be a smooth function. We use theorem \ref{Had} for $n=2$ and $p$ some point in $\mathbb{R}^k$. We then obtain the formula for $f\circ +:\mathbb{R}^k\times \mathbb{R}^k\rightarrow \mathbb{R}$:
 \[f(p+w)=f(p) +\sum_{l=1}^k w_l\cdot (\partial f/\partial x_l)(p) + \sum_{(i,j)\in \{1,...,k\}^2} w_i\cdot w_j\cdot h_{ij}(p,w).
\] for smooth functions $h_{ij}$. Thus, if $r=(r_i\oplus 0; 1\leq i\leq k)$ and $m=(0\oplus m_i; 1\leq i\leq k)$ then $M(f\circ +)(r,m)=$
\[M(f)(r)+\sum_{l=1}^k m_l\cdot M(\partial f/\partial x_l)(r) + \sum_{(i,j)\in k\times k} m_i\cdot m_j\cdot M(h_{ij})(y,w).\] The last term in this formula is zero, by definition of the ring structure. Since $M$ is a $C^\infty$-ring over $R$ we may replace the occurences of $M$ in other terms on the right hand side by $R$ to obtain:
\begin{equation}\label{clfrm}M(f)(r_i \oplus  m_i; 1\leq i\leq k)=R(f)(r)\oplus \sum_{l=1}^k m_l\cdot R(\partial f/\partial x_l)(r).\end{equation} This gives the closed formula I promised at the start of the proof and injectivity on objects as a corollary.\\
For essential surjectivity, just check that the above closed formula for $M$ yields a $C^\infty$-ring for each module $N$ over $U(R)$; the abelian group structure is trivially present and the underlying module is isomorphic to $N$.\\
The functor $\mathcal{U}$ is obviously faithful and therefore injective on arrows. To show that $\mathcal{U}$ is full, consider, for $\phi:M\rightarrow N$ a $U(R)$-linear map the assignment $\mathbb{R}^n\mapsto \phi^n$. Once it is checked that it is natural using the above closed formula one immediately sees that it is sent by $\mathcal{U}$ to $\phi$. $\Box$}

Recall that one may construct the pullback of a diagram $\xymatrix@1{& A\ar[d]^f \\ B \ar[r]_g & C}$ in categories as the subcategory of $A\times B$ containing only those arrows and objects $x$ for which $f(x)=g(x)$.

\bp{\textbf{Theorem.}\label{mainsex} $T_{C^\infty\text{-ring}}\simeq C^\infty\text{-ring}\times_{\text{CRing}}T_\text{CRing}$.\\
\underline{Proof:} On objects an inverse is provided by proposition \ref{Coomods}. As for arrows, let $(\alpha, (f,\phi)):(R,\mathfrak{A})\rightarrow (S,\mathfrak{B})$ be an arrow in the pullback. We seek an arrow $(f', \phi'):\mathcal{U}_R^{-1}(\mathfrak{A})\rightarrow \mathcal{U}_S^{-1}(\mathfrak{B})$. Take $f':=\alpha$ and define $\phi'_\mathbb{R}$ to be $\phi$. We must check a naturality square for $g:\mathbb{R}^n\rightarrow \mathbb{R}$ an arbitrary smooth function. This is easy using the closed formula (\ref{clfrm}) in the proof of proposition \ref{Coomods}:
\begin{multline}\notag \phi\left(R(g)(r)\oplus \sum_{l=1}^k m_l \cdot R(\partial f/\partial x_l)(r)\right)=
R(g)(r)\oplus \sum_{l=1}^k \phi(m_l)\cdot R(\partial f/\partial x_l)(r),\end{multline} as desired. $\Box$}

\bp{\textbf{Corollary.}\label{impcor} For any (co-)simplicial $C^\infty$-ring $R$, the category of abelian groups over $R$ is isomorphic to the category of (co-)simplicial modules of the underlying (co-)simplicial ring. $\Box$}

Let me now give some additional motivation for the notion of module given here. Call an abelian group object over some $C^\infty$-ring a module. Then derivations, after Quillen \cite{Q}, are sections of the projection onto the $C^\infty$-ring. One easily computes that for the case of commutative rings this yields derivations. In the case of $C^\infty$-rings we recover Dubuc and Kock's notion of module in \cite{DK}. In that text a module is just a module over the underlying commutative ring. They also define a derivation to be a linear map $d:R\rightarrow M$ for some $R$-module $M$ such that for every smooth function $f:\mathbb{R}^n\rightarrow \mathbb{R}$ we have that 
\[d(R(f)(r_1,...,r_n))=\sum_{i=1}^n R\left(\frac{\partial f}{\partial x_i}\right)\left(r_1,...,r_n\right) \cdot d(r_i).\]
This, they show, coincides with Quillen's notion of derivation applied to $C^\infty$-rings.

\end{document}